\documentclass[11pt]{article}
\usepackage{amssymb}
\usepackage{amsmath}
\usepackage{graphicx}
\usepackage{acronym}
\usepackage{amsfonts}
 \usepackage{cite}
 \usepackage{epstopdf}

\newtheorem{theorem}{Theorem}
\newtheorem{definition}{Definition}
\newtheorem{example}{Example}
\def\be{\begin{example}}
\def\ee{\end{example}}
\def\bt{\begin{theorem}}
\def\et{\end{theorem}\bigskip}
\def\bl{\begin{Lemma}}
\def\el{\end{Lemma}\bigskip}
\def\ep{\end{Proposition}\bigskip}
\def\bp{\begin{Proposition}}
\def\bd{\begin{definition}}
\def\ed{\end{definition}}

\newcommand{\alglist}{
\begin{list}{Step 1}
{\setlength{\leftmargin}{1.1 in}\setlength{\labelwidth}{1.0 in}} }

\newcommand{\x}{{\bf x}}

\newcommand{\uu}{{\bf u}}
\newcommand{\vv}{{\bf v}}

\newcommand{\A}{{\cal A}}

\begin{document}
\title{\bf  A Note on The Multidimensional Moment Problem}
 \author{Liqun Qi\footnote{Department of Applied
 Mathematics, The Hong Kong Polytechnic University, Hung Hom, Kowloon, Hong Kong.
His work was supported by the Hong Kong Research Grant Council
(Grant No. PolyU 501212, 501913, 15302114 and 15300715). E-mail:
maqilq@polyu.edu.hk.}}

\date{\today}
\maketitle

{\large

\begin{abstract}
In this note, we show that if a multidimensional sequence generates Hankel tensors and all the Hankel matrices, generated by this sequence, are positive semi-definite, then this sequence is a multidimensional moment sequence.
\medskip


\end{abstract}

\medskip

\section{Introduction}

The multidimensional moment problem is an important topic in mathematics \cite{Be, BCJ, Ham, Hav, Hav1, KM, Re, Va}.  In this note, we show that if a multidimensional sequence generates Hankel tensors and all the Hankel matrices, generated by this sequence, are positive semi-definite, then this sequence is a multidimensional moment sequence.   We do this in Section 2.   Some further questions are raised in Section 3.

We use small letters for scalars, bold letters for vectors, capital letters for matrices, and calligraphic letters for tensors.

\section{The Multidimensional Moment Problem}\label{sec:moment}

Denote $Z$ for the set of all positive integers, and $Z_+$ as the set of all nonnegative integers.   According to \cite{BCJ, Re}, a multidimensional sequence
\begin{equation} \label{e1}
S = \{ b_{j_1\ldots j_{n-1}} : j_1, \ldots, j_{n-1} \in Z_+ \}
\end{equation}
is called a {\bf multidimensional moment sequence} if there is a nonnegative measure $\mu$ on $\Re^{n-1}$ satisfying:
\begin{equation} \label{e2}
b_{j_1\ldots j_{n-1}} = \int_{\Re^{n-1}} t_1^{j_1}\ldots t_{n-1}^{j_{n-1}},  \ {\rm for}\  j_1, \ldots, j_{n-1} \in Z_+,
\end{equation}
are all finite.  For a given multidimensional sequence $S$ defined by (\ref{e1}), is it a multidimensional moment sequence? i.e., Is there a nonnegative measure $\mu$ such that (\ref{e2}) holds?  This problem is called the  {\bf multidimensional moment problem} \cite{Be, BCJ, Re}.

For any $m \in Z$, we may define a homogeneous polynomial of $n$ variables and degree $m$:
\begin{eqnarray} 
f(\x) & = & \sum \{ b_{j_1\ldots j_{n-1}}{m! \over j_1!\ldots j_{n-1}! (m-j_1-\ldots -j_{n-1})!} x_1^{j_1}\ldots x_{n-1}^{j_{n-1}}x_0^{m-j_1-\ldots - j_{n-1}}  \nonumber \\
&& :  j_1, \ldots, j_{n-1} \ge 0, j_1 + \ldots + j_{n-1} \le m \}. \label{e3}
\end{eqnarray}
According to \cite{Re}, $S$ is a multidimensional moment sequence if and only if for all $m$, $f(\x)$ a {\bf sum of $m$th power (SOM)} form.

A homogeneous polynomial $f(\x)$ of $n$ variables and degree $m$ is corresponding to an $m$th order $n$-dimensional symmetric tensor $\A = (a_{i_1\ldots i_m})$, where
\begin{equation} \label{e4}
a_{i_1\ldots i_m} = b_{j_1\ldots j_{n-1}},
\end{equation}
for  $j_{n-1} \ge 0, j_1 + \ldots + j_{n-1} \le m$, if in $\{ i_1, \cdots, i_m \}$, the frequency of $k$ is exactly $j_k$, $k= 1, \cdots, n-1$.
Then $f(\x)$ is an SOM form if and only if there are vectors $\uu_k \in \Re^n$ for $k = 1 \ldots, r$ such that
\begin{equation} \label{e5}
\A = \sum_{k=1}^r \uu_k^m,
\end{equation}
where for a vector $\vv \in \Re^n$, $\vv^m = (v_{i_1}\ldots v_{i_m})$ denotes a symmetric rank-one tensor.  Such a symmetric tensor is called a {\bf completely decomposable tensor} in \cite{LQX}.

Thus, a given multidimensional sequence $S$ defined by (\ref{e1}), is a multidimensional moment sequence if and only if all the symmetric tensors $\A$ generated by it are completely decomposable tensors for all $m$.   Note that when $m$ is odd, a symmetric tensor is always completely decomposable \cite{LQX}.

Suppose now that for $j_1, \ldots j_{n-1}, l_1, \ldots, l_{n-1} \in Z_+$, we have
\begin{equation} \label{e6}
b_{j_1\ldots j_{n-1}} = b_{l_1\ldots l_{n-1}}
\end{equation}
if 
\begin{equation} \label{e6a}
j_1 + 2j_2 + \ldots + (n-1)j_{n-1} =  l_1 + 2l_2 + \ldots + (n-1)l_{n-1}.
\end{equation}

By (\ref{e4}), for $i_1, \ldots i_n, k_1, \ldots, k_n \in Z_+$, we have
\begin{equation} \label{e7}
a_{i_1\ldots i_m} = a_{k_1\ldots k_m}
\end{equation}
as long as
\begin{equation} \label{e8}
i_1 + \ldots + i_m = k_1 + \ldots + k_m.
\end{equation}
By \cite{CQW, CQW1, DQW, LQW, LQX, Qi}, such a tensor is called a {\bf Hankel tensor}.   Thus, we call a multidimensional sequence $S$ satisfying (\ref{e8}) a {\bf Hankel multidimensional sequence}.

For an $m$th order $n$-dimensional Hankel tensor $\A = (a_{i_1\ldots i_m})$, by \cite{Qi}, there is a {\bf generating vector} $\vv = (v_0, \ldots, v_{mn})^\top$ such that
\begin{equation} \label{e9}
a_{i_1\ldots i_m} = v_{i_1+ \ldots + i_m}.
\end{equation}
If $m$ is even, then $\vv$ also generates a Hankel matrix $A$.  If the associated Hankel matrix $A$ is positive semi-definite, such a Hankel tensor $\A$ is called a {\bf strong Hankel tensor}.   By \cite{LQX}, a strong Hankel tensor is completely decomposable.   An explicit decomposition expression of a strong Hankel tensor is given in \cite{DQW}.

Furthermore, by (\ref{e6}), we see that 
\begin{equation} \label{e10}
v_{j_1+2j_2+\ldots +(n-1)j_{n-1}} =b_{j_1\ldots j_{n-1}}, 
\end{equation}
for $j_1, \ldots, j_{n-1} \in Z_+$, i.e., the components of $\vv$ are independent from $m$.   Thus, (\ref{e10}) defines an infinite sequence $V = \{ v_k : k \in Z_+ \}$.   This infinite sequence 
$V$ generates a sequence of Hankel matrices $H_p = (h_{ij})$, with $i, j = 0, \ldots, p-1, p \in Z$, and 
\begin{equation} \label{e11}
h_{ij} = v_{i+j}
\end{equation}\
for $i, j \in Z_+$.

By these, we have the following theorem.

\begin{theorem} \label{t1}
Suppose that a given multidimensional sequence $S$ defined by (\ref{e1}), satisfies (\ref{e6}), i.e., it is a Hankel  multidimensional sequence.   If all the Hankel tensors generated by $V$ are strong Hankel tensors, i.e., all the Hankel matrices $H_p$ generated by the sequence $V$ are positive semi-definite are positive semi-definite, then $S$ is a multidimensional moment sequence.
\end{theorem}

This links the classical result for the Hamburger moment problem \cite{Ham, Re}, and gives an application of the results in \cite{DQW, LQX, Qi}.

We may call such a sequence $V$ a {\bf strong Hankel sequence}.   The well-known strong Hankel sequence is the Hilbert sequence $\{ 1, {1 \over 2}, {1 \over 3}, \ldots \}$ \cite{SQ}.  

\section{Some Further Questions}\label{sec:dis}

The first question is: Are there any other Hankel multidimensional moment sequences, for which that infinite dimensional Hankel matrix $H$ is not positive semi-definite?    Another way to ask this question is as follows: Is there an infinite sequence $V = \{ v_k : k \in Z_+ \}$, not all the Hankel matrices $H_p$ generated by it are positive semi-definite, but all the $n$-dimensional Hankel tensors generated by it, for any order $m \in Z$, are completely decomposable?  Here $n \ge 3$ is a fixed positive integer.   Thus, the simplest case is that $n=3$.    We thus may ask:  Is there an infinite sequence $V = \{ v_k : k \in Z_+ \}$, not all the Hankel matrices $H_p$ generated by it are positive semi-definite, but all the $3$-dimensional Hankel tensors, generated by it, for any order $m \in Z$ are completely decomposable?   In \cite{LQW}, a class of sixth order three dimensional truncated Hankel tensors are discussed.   Such Hankel tensors are not strong Hankel tensors \cite{LQW, LQX}, but still completely decomposable \cite{LQX}.   Can we build an infinite sequence $V$ to answer this question, based upon such sixth order three dimensional truncated Hankel tensors?   

If such a sequence $V$ exists, then the next question is: what are the necessary and sufficient conditions such a sequence $V$ should satisfy?   

Note that until now, all the known positive semi-definite Hankel tensors are SOS \cite{CQW, LQW, LQX}, but there are positive semi-definite Hankel tensors which are not completely decomposable.  This situation may make such a characterization somewhat subtle.  

Another question is: Are there multidimensional moment sequences, which are not Hankel multidimensional sequences, i.e., condition (\ref{e6}) is not satisfied.   If the answer to this question is ``yes'', then how to characterize such multidimensional sequences?

{\bf Acknowledgment}  The author is thankful to Weiyang Ding for his comments.

 }


\end{document}